\documentclass[12pt]{amsart}
\usepackage[babel]{csquotes}
\usepackage{enumitem}
\usepackage{amsmath,amsthm,amssymb,mathrsfs,amsfonts,verbatim,enumitem,color,leftidx}
\usepackage{kotex}
\usepackage{mathabx}
\usepackage{mathtools}
\usepackage{graphicx}
\usepackage[left=2.9cm,right=2.9cm,top=3.3cm,bottom=3.4cm,a4paper]{geometry}
\usepackage{etoolbox} 
\usepackage{tikz}
\usepackage{bbm}
\usepackage[all,tips]{xy}
\usepackage{graphicx,ifpdf}
\usepackage{stmaryrd}
\ifpdf
\fi

\usepackage[right,displaymath,mathlines]{lineno}
\allowdisplaybreaks[0]

\definecolor{dartmouthgreen}{rgb}{0.05, 0.5, 0.06}

\usepackage[colorlinks]{hyperref}
\hypersetup{
linkcolor=blue,          
citecolor=dartmouthgreen,      
}

\usepackage{tikz}
\usepackage{here}
\usetikzlibrary{arrows,snakes,backgrounds}
\usetikzlibrary{arrows.meta,calc}

\newtheorem{thm}{Theorem}[section]
\newtheorem{lem}[thm]{Lemma}
\newtheorem{coro}[thm]{Corollary}
\newtheorem{prop}[thm]{Proposition}

\theoremstyle{definition}

\newtheorem{remark}[thm]{Remark}

\theoremstyle{remark}

\numberwithin{equation}{section}

\def\wt{\widetilde}

\def\mc{\mathcal}

\def\G{\mathrm{G}}
\def\H{\mathrm{H}}

\def\A{\mathbb{A}}

\def\e{\epsilon}

\def\CC{\mathbb{C}}

\def\Ind{\operatorname{Ind}}

\def\Res{\operatorname{Res}}

\def\beq{\begin{equation}}
\def\eeq{\end{equation}}
\def\beqn{\begin{equation*}}
\def\eeqn{\end{equation*}}
\def\beqna{\begin{eqnarray}}
\def\eeqna{\end{eqnarray}}
\def\beqnan{\begin{eqnarray*}}
\def\eeqnan{\end{eqnarray*}}

\definecolor{esperance}{rgb}{0.0,0.5,0.0}
\definecolor{darkspringgreen}{rgb}{0.09, 0.45, 0.27}
\definecolor{princetonorange}{rgb}{1.0, 0.56, 0.0}

\title[Transition probability]{Transition probability of discrete geodesic flow on the standard non-uniform quotient of $PGL_3$}

\begin{document}

\begin{abstract}
We describe the local transition probability of a singular diagonal action on the standard non-uniform quotient of $PGL_3$ associated to the type 1 geodesic flow. As a consequence, we deduce the property of strongly positive recurrence.
\end{abstract}


\author{Sanghoon Kwon}
\address{Sanghoon Kwon*}
\curraddr{Room 506 Department of Mathematical Education\\ Catholic Kwandong University \\ Gangneung 25601 \\ Republic of Korea}
\email{skw1.math@gmail.com, skwon@cku.ac.kr}

\thanks{This work has been supported by NRF grant (No. RS-2023-00237811).}


\maketitle
\tableofcontents

\section{Introduction}

In this paper, we discuss the strongly positive recurrence property of discrete geodesic flows in the standard arithmetic quotient of the affine building for $PGL_3$. Strongly positive recurrence property of countable topological Markov chains was defined in \cite{Ru}. If a directed graph of a Markov chain is not strongly positive recurrent, then the entropy is mainly concentrated near infinity in the sense that it is supported by the infinite paths that spend most of their time outside a finite subgraph. Recently, \cite{RS} proved the effective intrinsic ergodicity for all strongly positive recurrent topological Markov shifts. Namely, they provide an effective bound of the distance between an invariant measure and the measure of maximal entropy in terms of the difference of their entropies. In \cite{GST}, the authors investigate the notion of strongly positive recurrence of geodesic flows on non-compact negatively curved manifolds, using the entropy and pressure at infinity. Our setting, the geodesic flow on an affine building of rank 2, can be viewed as an example of a space with non-positive curvature.

Let $\mathbb{F}_q$ be the finite field of order $q$ and let $\mathbb{F}_q(t)$ be the field of rational functions over $\mathbb{F}_q$. The absolute value $\|\cdot\|$ of $\mathbb{F}_q(t)$ is defined for any $f\in \mathbb{F}_q(t)$, by 
$$\|{f}\|:=q^{\deg (g)-\deg (h)},$$
for $g,h$ are polynomial over $\mathbb{F}_q$ satisfying $f=\frac{g}{h}$.
The completion of $\mathbb{F}_q(t)$ with respect to $\|\cdot\|$, the field of formal Laurent series in $t^{-1}$, is denoted by $\mathbb{F}_q(\!(t^{-1})\!)$, i.e.,
$$\mathbb{F}_q(\!(t^{-1})\!):=\left\{\sum_{n=-N}^\infty a_nt^{-n}:N\in \mathbb{Z}, a_n\in \mathbb{F}_q\right\}.$$ 
The valuation ring $\mathcal{O}$ is the subring of power series 
$$\mathbb{F}_q\mathbb{[\![}t^{-1}]\!]:=\left\{\sum_{n=0}^\infty a_nt^{-n}: a_n\in \mathbb{F}_q\right\}.$$

Let $G$ be the group $PGL(3,\mathbb{F}_q(\!(t^{-1})\!))$, $\Gamma$ be the standard non-uniform arithmetic lattice $PGL(3,\mathbb{F}_q[t])$ of $G$ and $K$ be a maximal compact subgroup $PGL(3,\mathcal{O})$ of $G$. Denote by $\mathcal{B}$ the building $\mathcal{B}_3(\mathbb{F}_q(\!(t^{-1})\!))$ associated to the group $G$. It is the 2-dimensional contractible simplicial complex defined as follows. We say two $\mathcal{O}$-lattices $L$ and $L'$ of rank $3$ are in the same equivalence class if $L=sL'$ for some $s\in \mathbb{F}_q(\!(t^{-1})\!)^\times$. The set of the equivalence classes $[L]$ will be the set of vertices of $\mathcal{B}$. For given $k$-vertices $[L_1],[L_2],[L_3]$, they form a $2$-dimensional simplex in $\mathcal{B}$ if 
\begin{equation*}\label{eq:1.1}
t^{-1} \Lambda_1\subset \Lambda_3\subset \Lambda_2\subset \Lambda_1
\end{equation*}
for some $\Lambda_i\in [L_i].$ Then, the set of vertices of $\mathcal{B}$ may be identified with $G/K$.  

For a more comprehensive and detailed discussion about Bruhat-Tits building, one may follow \cite{AB}. We also remark that \cite{Mo} investigated the dynamical properties of the diagonal action in the compact quotient of a $p$-adic Chevalley group. This paper explores one of the non-compact generalization of the dynamical system discussed in \cite{Mo}.

The \emph{type} $\tau(x)$ of vertex $x=gK$ is defined by $\log_q\!\|\det(g)\|\textrm{ (mod }3$). Each apartment of $\mathcal{B}$ is a Euclidean plane tiled with equilateral triangles. The type $\tau(v\to w)$ of a directed edge $v\to w$ from a vertex $v$ to a vertex $w$ is defined to be $\tau(w)-\tau(v)$. If $e$ is a directed edge from $v$ to $w$ in $\mathcal{B}$, then we denote by $s(e)=v$ (source) and $t(e)=w$ (target). A sequence of $e_1,e_2,\ldots,e_n$ of directed edges in $\mathcal{B}$ is called a path if $t(e_k)=s(e_{k+1})$ for all $1\le k\le n-1$. If it consists of type $i$ directed edges, then it is called a path of type $i$. A path $e_1,e_2,\ldots,e_n$ in $\mathcal{B}$ is called a \emph{geodesic path} if it is a part of straight line in an apartment in $\mathcal{B}$. Equivalently, it is a path with the condition that $s(e_k),s(e_{k+1})=t(e_k),t(e_{k+1})$ do not form a chamber in $\mathcal{B}$ for all $1\le k\le n-1$. See Figure~\ref{fig:geopath}.

\begin{center}
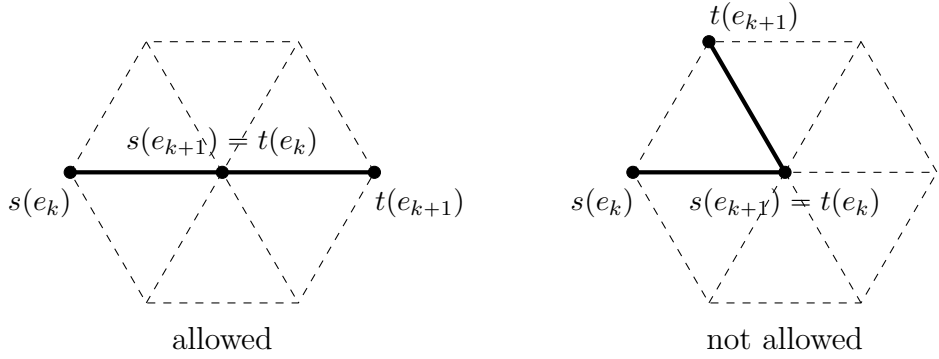
\begin{figure}[h]
\label{fig:geopath}
\newcommand*\rows{6}
\begin{tikzpicture}[scale=2]

\draw[dashed] (0:1) --++(120:1) --++(180:1) --++(240:1) --++ (300:1) --++(0:1) --++ (60:1);
\draw[dashed] (180:1) --++ (0:2);
\draw[dashed] (240:1) --++ (60:2);
\draw[dashed] (120:1) --++ (300:2);

\draw[ultra thick] (180:1) --++ (0:2);
\draw[ultra thick] (2.7,0) --++ (0:1) --++ (120:1);

\draw[dashed] (2.7,0) --++(0:2) --++(120:1) --++ (180:1) --++(240:1) --++ (300:1) --++(0:1) --++ (60:1);
\draw[dashed] (3.7,0) --++ (60:1);
\draw[dashed] (3.7,0) --++ (120:1);
\draw[dashed] (3.7,0) --++ (240:1);
\draw[dashed] (3.7,0) --++ (300:1);

\draw[fill] (0,0) circle (0.04cm);
\draw[fill] (-1,0) circle (0.04cm);
\draw[fill] (1,0) circle (0.04cm);
\draw[fill] (2.7,0) circle (0.04cm);
\draw[fill] (3.7,0) circle (0.04cm);
\draw[fill] (3.2,{sqrt(3)/2}) circle (0.04cm);

\node (a) at (270:1.1) {allowed};
\node (b) at (3.7,-1.1) {not allowed};
\node (c) at (-1.2,-0.2) {\small$s(e_k)$};
\node (d) at (0,0.2) {\small$s(e_{k+1})=t(e_k)$};
\node (e) at (1.3,-0.2) {\small$t(e_{k+1})$};
\node (f) at (2.5,-0.2) {\small$s(e_k)$};
\node (g) at (3.7,-0.2) {\small$s(e_{k+1})=t(e_k)$};
\node (h) at (3.5,1.03) {\small$t(e_{k+1})$};
\end{tikzpicture}
\caption{admissible geodesic paths}\label{fig:geopath}
\end{figure}
\end{center}

As we mentioned earlier, we explore the recurrence property of the type $1$ geodesic flow in $\Gamma\backslash \mathcal{B}$, that is, the shift map $[(e_n)_{n\in\mathbb{Z}}]_\Gamma\mapsto [(e_{n+1})_{n\in\mathbb{Z}}]_\Gamma$ for type 1 geodesics $(e_n)_{n\in\mathbb{Z}}$ in $\mathcal{B}$.
Consider a standard type $1$ bi-infinite geodesic $\mathbf{s}=(\ldots,v_{-1},v_0,v_{1},\ldots)$ where $v_n=\textrm{diag}(t^n,1,1)K$ in $\mathcal{B}$. We observe that
\[g \begin{pmatrix}t^n & 0 & 0 \\ 0 & 1 & 0 \\ 0 & 0 & 1\end{pmatrix}K=\begin{pmatrix}t^n & 0 & 0 \\ 0 & 1 & 0 \\ 0 & 0 & 1\end{pmatrix}K\]
for all $n\in\mathbb{Z}$ if and only if $g$ is an element of $K$ of the form \[g=\begin{pmatrix}* & 0 & 0 \\ 0 & * & * \\ 0 & * & *\end{pmatrix}.\] 
Thus, if we denote by $a$ the element $\textrm{diag}(t,1,1)$ in $G$, then the type $1$ geodesic flow system corresponds to the right multiplication action $T_a\colon \Gamma\backslash G/M\to \Gamma\backslash G/M$ given by $T_a(\Gamma gM)=\Gamma gaM$ for \[M=\left\{m=\begin{pmatrix}k_{11} & 0 & 0 \\ 0 & k_{22} & k_{23} \\ 0 & k_{32} & k_{33} \end{pmatrix}\colon m\in K\right\}.\]  We will investigate the asymptotic behavior of the number of periodic orbits of the system $T_a\colon \Gamma\backslash G/M\to\Gamma\backslash G/M$.

Let $\pi_K\colon \Gamma\backslash G/M\to \Gamma\backslash G/K$ be the natural projection map and denote by $o$ the identity coset $\Gamma eK$ in $\Gamma\backslash G/K$. Let $f_n(o)$ denote the number of first return cycles at $o$ of length $n$. Namely,
 \[f_n(o)=\#\big\{\mathbf{x}\in \Gamma\backslash G/M\colon \pi_K(\mathbf{x})=o,T_a^n(\mathbf{x})=\mathbf{x},\, n=\min \{k>0\colon \pi_K(T_a^k(\mathbf{x}))=o\}\big\}.\]
Also for each $x\in \Gamma\backslash G/K$, let $g_n(x)$ denote the number 
 \[g_n(x)=\#\{\mathbf{x}\in \Gamma\backslash G/M\colon \pi_K(\mathbf{x})=x,\, T_a^n(\mathbf{x})=\mathbf{x}\}\]
of closed cycles based at $x$ of length $n$.
 
Additionally, since the period of $T_a$ on $\Gamma\backslash G/M$ is 3, the Gurevich entropy $h$ (based at $x$) of $T_a$ is defined by
 \[h_{T_a}=\lim_{n\to\infty}\frac{1}{3n}\log g_{3n}(x).\]
This value does not depend on the choice of $x$. The following is the main theorem of this article.

\begin{thm}\label{thm:main}
The type 1 discrete geodesic flow system $(\Gamma\backslash G/M,T_a)$ is strongly positive recurrent in the sense that
\[\limsup_{n\to\infty}\frac{1}{n}\log f_n(o)<h_{T_a}.\]
\end{thm}
We prove this theorem in Section~\ref{sec:3}. In fact, we prove that $h_{T_a}=2\log q$ and $\underset{n\to\infty}{\limsup}\,\frac{1}{n}\log f_n(o)=\frac{5}{3}\log q$.
\bigskip

It is worth mentioning that \cite{AGP} demonstrate the logarithmic law of geodesic flow in the non-compact quotient of affine buildings. It would be an interesting result if we could get a theorem on the limiting distribution of extreme values, which is obtained in \cite{KL1} in the case of geometrically finite quotients of trees.
\bigskip

\section{Reduction to countable Markov Shift}\label{sec:2}

Recall that $\mathbf{s}=(\ldots,v_{-1},v_0,v_{1},\ldots)$ is defined by the sequence of vertices $v_n=\textrm{diag}(t^n,1,1)K$ of the \emph{standard} type 1 geodesic in $\mathcal{B}$. Let $(Y,\sigma)$ be the shift space given by 
\[Y=\{\mathbf{y}\in (G/K)^{\mathbb{Z}}\colon \mathbf{y}=(y_n)_n\textrm{ corresponds to a type 1 geodesic in }\mathcal{B}\}\]
and $\sigma\colon Y\to Y$, $\sigma((y_n)_n)=(y_{n+1})_n$.
Let $\Phi\colon G/M\to Y$ be the bijective map given by $\Phi(gM)=g\mathbf{s}=(\ldots,gv_{-2},gv_{-1},gv_{0},gv_{1},gv_{2},\ldots)$. Then, we have $\Phi\circ T_a=\sigma\circ\Phi$.
\begin{equation*}
\begin{aligned}
\xymatrix{
G/M \ar[r]^{\,\,T_a}\ar[d]_{\Phi} & G/M\ar[d]_{\Phi}\\
Y \ar[r]^{\,\,\sigma} & Y }
\end{aligned}
\end{equation*}

Now let $X=Y/\sim$ where $\mathbf{y}\sim\mathbf{y}'\Leftrightarrow \mathbf{y}=\gamma\mathbf{y}'$ for some $\gamma\in \Gamma$. Then, the following diagram also commutes.
\begin{equation*}
\begin{aligned}
\xymatrix{
\Gamma\backslash G/M \ar[r]^{T_a}\ar[d]_{\phi} & \Gamma\backslash G/M \ar[d]_{\phi}\\
X\ar[r]^{\sigma}& X }
\end{aligned}
\end{equation*}

Let 
\[I=\left\{\begin{pmatrix}k_{11} & k_{12} & k_{13} \\ k_{21} & k_{22} & k_{23} \\ k_{31} & k_{32} & k_{33}\end{pmatrix}\in K\colon \|k_{21}\|,\|k_{31}\|\le q^{-1}\right\}\] be the Iwahori subgroup of $G$, which is the stabilizer of the type 1 directed edge $K\to aK$ in $\mathcal{B}$. Then, the set of type 1 directed edges can be identified with $G/I$.

Now let $\mathcal{E}=\{\mathbf{e}\in(G/I)^\mathbb{Z}\colon \mathbf{e}=(e_n)_n\textrm{ corresponds to a type 1 geodesic in }\mathcal{B}\}$. In other words, $(e_n)_n\in\mathcal{E}$ if $t(e_n)=s(e_{n+1})$ and $s(e_n)$, $s(e_{n+1})$, $t(e_{n+1})$ do not form a chamber in $\mathcal{B}$ for all $n\in\mathbb{Z}$. Then, the map $\Psi\colon Y\to \mathcal{E}$ given by 
\[\Psi((y_n)_{n\in\mathbb{Z}})=(e_n)_{n\in\mathbb{Z}} \quad\textrm{ (where }s(e_n)=y_n,t(e_n)=y_{n+1}\textrm{)}\] will be the natural bijection. Now, we say that a sequence $\mathbf{d}\in (\Gamma\backslash G/I)^{\mathbb{Z}}$ of directed edges in $\Gamma\backslash\mathcal{B}$ is \emph{admissible} if $\mathbf{d}=(d_n)_n$ may lift to an element in $\mathcal{E}$. Denote by $\mathcal{D}$ the set $\{\mathbf{d}\in (\Gamma\backslash G/I)^{\mathbb{Z}}\colon \mathbf{d}\textrm{ is admissible}\}$ of admissible sequences. If we define the equivalence relation on $\mathcal{E}$ by $\mathbf{e}\sim\mathbf{e}'\Leftrightarrow \mathbf{e}=\gamma\mathbf{e}'$ for some $\gamma\in \Gamma$, then the map $\Psi$ induces a surjection $\psi$ from $X$ to $\mathcal{D}$. The relation $p_3\circ \Psi=\psi\circ p_2$ also holds.

Let us also denote by $p_1,p_2,p_3$ the projection map $G/M\to \Gamma\backslash G/M$, $Y\to X$, and $\mathcal{E}\to\mathcal{D}$, respectively. The following commutative diagram describes the notations.
\begin{equation*}
\begin{aligned}
\xymatrix{
G/M \ar[r]_{\,\,\,\Phi}^{\,\,\,\simeq}\ar[d]_{p_1} & Y\ar[r]_{\Psi}^{\simeq}\ar[d]_{p_2} &\mathcal{E}\ar[d]_{p_3}\\
\Gamma\backslash G/M\ar[d]_{\pi_K}\ar[r]_{\quad\phi}^{\quad\simeq} & X\ar@{.>}[r]_{\psi}^{} &\mathcal{D}  \\
\Gamma\backslash G/K &
}
\end{aligned}
\end{equation*}


We recall that
 \[g_n(o)=\#\{\mathbf{x}\in \Gamma\backslash G/M\colon \pi_K(\mathbf{x})=o,\, T_a^n(\mathbf{x})=\mathbf{x}\}\]
and 
\[f_n(o)=\#\big\{\mathbf{x}\in \Gamma\backslash G/M\colon \pi_K(\mathbf{x})=o,T_a^n(\mathbf{x})=\mathbf{x},\, n=\min \{k>0\colon \pi_K(T_a^k(\mathbf{x}))=o\}\big\}.\]
Let us define
\begin{align*}
\mathcal{D}_{\textrm{per},n}=&\,\{\mathbf{d}\in\mathcal{D}\colon d_0=\Gamma eI,\sigma^n(\mathbf{d})=\mathbf{d}\}\\
\mathcal{D}_{\textrm{prim},n}=&\,\big\{\mathbf{d}\in\mathcal{D}\colon d_0=\Gamma eI,\sigma^n(\mathbf{d})=\mathbf{d},\, n=\min\{k>0\colon d_k=\Gamma eI\}\big\}.
\end{align*}
\begin{lem}
We have 
\[g_n(o)=\sum_{\mathbf{d}\in\mathcal{D}_{\textrm{per},n}}\bigr|p_3^{-1}(\mathbf{d})\bigr|\qquad\textrm{and}\qquad f_n(o)=\sum_{\mathbf{d}\in\mathcal{D}_{\textrm{prim},n}}\bigr|p_3^{-1}(\mathbf{d})\bigr|.\]
\end{lem}
\begin{proof}
Using the isomorphism $\phi\colon \Gamma\backslash G/M\to X$, we may identify $\Gamma\backslash G/M$ to $X$. Let $\mathbf{x}=(x_n)_{n\in\mathbb{Z}}$ be an element in $X$ satisfying $\pi_K(\mathbf{x})=o$, $T_a^n(\mathbf{x})=\mathbf{x}$. It corresponds to a finite admissible sequence $(x_0,x_1,\ldots,x_{n-1},x_n)$ in $(G/K)^{n+1}/\sim$ for which $x_0\in G/K$ is a lift of $o= \Gamma eK$ and there exists $\gamma x_0=x_n$ for some $\gamma\in \Gamma$. Moreover, we may choose a unique representative $(v_0,y_1,\ldots,y_{n-1},\gamma v_0)$ in $(G/K)^{n+1}$. Note also that $\psi(\mathbf{x})$ is an element in $\mathcal{D}_{\textrm{per},n}$.

Conversely, let $\mathbf{d}\in \mathcal{D}_{\textrm{per},n}$. Since $\mathbf{d}$ is admissible, there is an $\mathbf{x}=(x_n)_{n\in\mathbb{Z}}\in X$ such that $\psi(\mathbf{x})=\mathbf{d}$. Furthermore, $\pi_K(\mathbf{x})=o$ and $T_a^n(\mathbf{x})=\mathbf{x}$. From the above observation, we have
\begin{align*}
&\big|\{\mathbf{x}\in X\colon\Psi(\mathbf{x})=\mathbf{d}\}\big|\\
=\,&\big|\{\overline{\mathbf{y}}=(v_0,y_1,\ldots,y_{n-1},\gamma v_0)\in (G/K)^{n+1}\colon \overline{\mathbf{y}}\textrm{ is admissible},\,p_3(\Psi(\overline{\mathbf{y}}))=(d_0,\ldots,d_{n})\}\big|\\
=\,&\big|\{\mathbf{y}\in Y\colon p_3(\Psi(\mathbf{y}))=\mathbf{d}\}\big|
\end{align*}
which yields 
\[\#\{\mathbf{x}\in \Gamma\backslash G/M\colon \pi_K(\mathbf{x})=o,\, T_a^n(\mathbf{x})=\mathbf{x}\}=\sum_{\mathbf{d}\in\mathcal{D}_{\textrm{per},n}}\bigr|p_3^{-1}(\mathbf{d})\bigr|.\]
Similar argument also gives the second equality.
\end{proof}

 \section{Proof of Theorem~\ref{thm:main}}\label{sec:3}
 
In this section, we prove Theorem~\ref{thm:main}. Namely, we show that \[\limsup_{n\to\infty}\frac{1}{n}\log f_n(o)<h_{T_a}=\lim_{n\to\infty}\frac{1}{3n}\log g_{3n}(o).\]
 The proof goes through the explicit calculation of $f_n(o)$ and $g_n(o)$ by investigating local transition probabilities of Markov shift $(\mathcal{D},\sigma)$.

First, the Birkhoff decomposition says that given every $g\in G$, there exists a unique pair of non-negative integers $(m,n)$ with $0\le n\le m$ such that
$$g\in\Gamma\begin{pmatrix}t^m & 0 & 0 \\ 0 & t^n & 0 \\ 0 & 0 & 1\end{pmatrix} K$$ holds. For the reduction algorithm, see Lemma~3.2 of \cite{HK}. Hence, we may denote by $v_{m,n}$ the vertex of the quotient complex $\Gamma\backslash \mathcal{B}(G)$ corresponds to
$$\Gamma\textrm{diag}(t^m,t^n,1)K.$$ There is an edge between two vertices $v_{m,n}$ and $v_{m',n'}$ if and only if the following hold:
\begin{displaymath}
\begin{cases} (m',n')\in\{(m\pm1,n),(m,n\pm1),(m\pm1,n\pm1)\} &\text{ if }m>n>0\\
(m',n')\in \{(m\pm1,n),(m,n+1),(m+1,n+1)\} &\text{ if }m>n=0\\
(m',n')\in \{(m+1,n),(m,n-1),(m\pm1,n\pm1)\} &\text{ if }m=n>0\\
(m',n')\in \{(1,0),(1,1)\}&\text{ if }m=n=0.
\end{cases}
\end{displaymath}
We denote by $e_{\frac{m+m'}{2},\frac{n+n'}{2}}$ the type $1$ directed edge from $v_{m,n}$ to $v_{m',n'}$. See Figure~\ref{fig:quotient} for the picture of the quotient complex $\Gamma\backslash G/K$.

\begin{center}
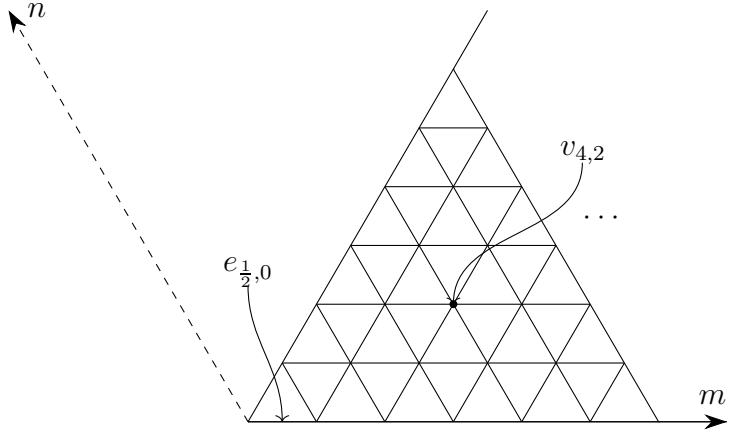
\begin{figure}[h]
\newcommand*\rows{6}
\begin{tikzpicture}[scale=0.9]
    \foreach \row in {0, 1, ...,\rows} {
        \draw ($\row*(0.5, {0.5*sqrt(3)})$) -- ($(\rows,0)+\row*(-0.5, {0.5*sqrt(3)})$);
        \draw ($\row*(1, 0)$) -- ($(\rows/2,{\rows/2*sqrt(3)})+\row*(0.5,{-0.5*sqrt(3)})$);
        \draw ($\row*(1, 0)$) -- ($(0,0)+\row*(0.5,{0.5*sqrt(3)})$);
    }
\draw (60:6) --++(60:1);
\draw (0:6) --++(0:1);
\draw[fill] (30:3.466) circle (0.05);
\node (a) at (30:6) {$\cdots$};
\node (b) at (3:6.8) {$m$};
\node (c) at (117:6.8) {$n$};
\node (d) at (39:6.3) {$v_{4,2}$};
\node (e) at (0,2.2) {$e_{\frac{1}{2},0}$};
\draw[->](38:6.2) to [out=270,in=90] node[above]{}(30:3.466);
\draw[->](0,2) to [out=270,in=90] node[above]{}(0:0.5);

\draw [-{Stealth[length=3mm, width=2mm]}] (120:6.9) -- (120:7);
\draw [dashed] (120:6.9) -- (0,0);
\draw [-{Stealth[length=3mm, width=2mm]}] (0,0) -- (0:7);

\end{tikzpicture}
\caption{$\Gamma\backslash G/K$}\label{fig:quotient}
\end{figure}
\end{center}

\begin{figure}[h]
\begin{center}
\begin{tikzpicture}[scale=0.97]
\fill (-4,0)[red] circle (2pt);\fill (0,0) circle (2pt);\fill [blue] (4,0) circle (2pt);\fill [blue] (-2,{sqrt(12)}) circle (2pt);\fill [red] (2,{sqrt(12)}) circle (2pt);\fill [blue] (-2,-{sqrt(12)}) circle (2pt);\fill [red] (2,-{sqrt(12)}) circle (2pt);\fill [red] (-2,-{sqrt(12)}/3) circle (2pt);\fill [blue] (-2,+{sqrt(12)}/3) circle (2pt);\fill [red] (2,+{sqrt(12)}/3) circle (2pt);\fill [blue] (2,-{sqrt(12)}/3) circle (2pt);\fill [red] (-1,+{sqrt(300)}/6) circle (2pt); \fill [blue] (0,+{sqrt(48)}/3) circle (2pt); \fill [red] (-1,-{sqrt(300)}/6) circle (2pt); \fill [blue] (2,-{sqrt(48)}/3) circle (2pt);
\draw [thick] (-4,0)--(-2,{sqrt(12)})--(2,{sqrt(12)})--(4,0)--(2,-{sqrt(12)})--(-2,-{sqrt(12)})--(-4,0);\draw [thick] (-4,0)--(-2,+{sqrt(12)}/3);\draw [thick] (-1,-{sqrt(300)}/6)--(0,0)--(0,+{sqrt(48)}/3)--(-1,-{sqrt(300)}/6); \draw [thick,shorten >=110pt] (-1,+{sqrt(300)}/6)--(2,-{sqrt(12)}/3); \draw [thick,shorten >=50pt] (2,-{sqrt(12)}/3)--(-1,+{sqrt(300)}/6); \draw [thick] (4,0)--(2,+{sqrt(12)}/3); \draw [thick] (-2,-{sqrt(12)}/3)--(-2,-{sqrt(12)});  \draw[thick,shorten >=16pt] (-2,+{sqrt(12)}/3)--(0,0); \draw [thick] (-2,+{sqrt(12)}/3)--(-1,-{sqrt(300)}/6);\draw [thick] (0,0)--(2,-{sqrt(48)}/3)--(-1,-{sqrt(300)}/6);\draw [thick] (-2,+{sqrt(12)})--(-1,+{sqrt(300)}/6); 
\draw[thick,shorten >=60pt] (2,-{sqrt(48)}/3)--(-1,+{sqrt(300)}/6);
\draw[thick,shorten >=134pt] (-1,+{sqrt(300)}/6)--(2,-{sqrt(48)}/3);\draw[thick,shorten >=82pt] (2,-{sqrt(12)})--(0,0);\draw [thick] (2,-{sqrt(12)})--(2,-{sqrt(48)}/3);\draw [thick] (2,+{sqrt(12)})--(0,+{sqrt(48)}/3); \draw[thick,shorten >=40pt] (-2,+{sqrt(12)})--(0,0);\draw[thick, shorten >=45pt] (-1,+{sqrt(300)}/6)--(0,0); \draw[thick, shorten >=36pt] (2,-{sqrt(12)}/3)--(0,0);
\draw[thick, shorten >=100pt] (-2,-{sqrt(12)}/3)--(2,-{sqrt(12)}/3);\draw[thick, shorten >=100pt] (-2,-{sqrt(12)}/3)--(2,-{sqrt(12)}/3);\draw[thick, shorten >=97pt] (2,-{sqrt(12)}/3)--(-2,-{sqrt(12)}/3); \draw[thick, shorten >=50] (-4,0)--(0,0);\draw[thick, shorten >=60] (4,0)--(0,0);\draw[thick](2,+{sqrt(12)}/3)--(2,-{sqrt(12)}/3);\draw[thick, shorten >=30] (2,+{sqrt(12)}/3)--(0,0);\draw[thick, shorten >=70] (-2,-{sqrt(12)})--(0,0);\draw [thick, shorten >=40] (2,+{sqrt(12)})--(0,0);\draw [thick,shorten >=93] (-2,-{sqrt(12)}/3)--(0,+{sqrt(48)}/3);\draw [thick,shorten >=60] (0,+{sqrt(48)}/3)--(-2,-{sqrt(12)}/3);\draw [thick,shorten >=67] (2,+{sqrt(12)}/3)--(-2,+{sqrt(12)}/3);\draw [thick,shorten >=80] (-2,+{sqrt(12)}/3)--(2,+{sqrt(12)}/3);\draw [thick,shorten >=50] (-2,-{sqrt(12)}/3)--(0,0);
\node at (2,1.4) {\tiny{$q-1$}};\node at (-2,1.4) {\tiny{$q-1$}};\node at (2.2,-1.4) {\tiny{$(q-1)^2$}};\node at (-2.4,-1.2) {\tiny{$q-1$}};\node at (-0.9,-3.1) {\tiny{$(q-1)^2$}};\node at (2.3,-2.1) {\tiny{$q-1$}};\node at (0.4,2.2) {\tiny{$q-1$}};\node at (-0.6,3) {\tiny{$q-1$}};

\draw[dashed] (-1,-{sqrt(299)}/6) -- (2,-{sqrt(12)}/3);
\end{tikzpicture}
\end{center}
\caption{Star of a vertex in $\mathcal{B}$}\label{fig:star}
\end{figure}
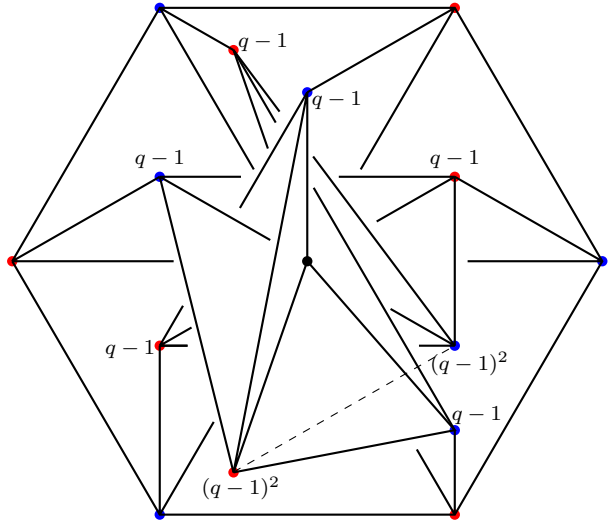

Meanwhile, for any vertex $x=gK\in G/K$ of type $i$ in $\mathcal{B}$, there are $q^2+q+1$ vertices of type $i+1$ neighbors of $gK$. They are given by
\[\left\{g\begin{pmatrix}t & 0 & 0 \\ 0 & 1 & 0 \\ 0 & 0 & 1\end{pmatrix} K\right\}\cup\left\{g\begin{pmatrix}1 & bt & 0 \\ 0 & t & 0 \\ 0 & 0 & 1\end{pmatrix} K\colon b\in\mathbb{F}_q\right\}\cup\left\{\begin{pmatrix}t^{-1} & 0 & c \\ 0 & t^{-1} & d \\ 0 & 0 & 1\end{pmatrix} K\colon c,d\in\mathbb{F}_q\right\}.\]
Also, there are $q^2+q+1$ vertices of type $i+2$ neighbors of $gK$, which are 
\[\left\{g\begin{pmatrix}t & 0 & 0 \\ 0 & t & 0 \\ 0 & 0 & 1\end{pmatrix} K\right\}\cup\left\{g\begin{pmatrix}1 & 0 & 0 \\ 0 & t^{-1} & b \\ 0 & 0 & 1\end{pmatrix} K\colon b\in\mathbb{F}_q\right\}\cup\left\{\begin{pmatrix}t^{-1} & c & d \\ 0 & 1 & 0 \\ 0 & 0 & 1\end{pmatrix} K\colon c,d\in\mathbb{F}_q\right\}.\]
See Figure~\ref{fig:star} for the star of a vertex in $\mathcal{B}$. A similar discussion with the picture is presented in Subsection~3.2 of \cite{GP} and Section~3 of \cite{KL2}.

\bigskip

\begin{figure}[h]
\begin{center}
\begin{tikzpicture}[scale=0.9]
\draw[-{Stealth[open, length=3mm, width=2mm]}] (-8,0)--(-6,0);
\draw[-{Stealth[length=3mm, width=2mm]}] (-6,0)--(-4,0);
\draw[-{Stealth[length=3mm, width=2mm]}] (-6,0)--(-7,{sqrt(3)});
\draw[-{Stealth[length=3mm, width=2mm]}] (-6,0)--(-7,-{sqrt(3)});
\draw[-{Stealth[open, length=3mm, width=2mm]}] (1,-{sqrt(3)})--(-0,0);
\draw[-{Stealth[length=3mm, width=2mm]}] (0,0)--(2,0);
\draw[-{Stealth[length=3mm, width=2mm]}] (0,0)--(-1,{sqrt(3)});
\draw[-{Stealth[length=3mm, width=2mm]}] (0,0)--(-1,-{sqrt(3)});
\draw[-{Stealth[open, length=3mm, width=2mm]}] (7,{sqrt(3)})--(6,0);
\draw[-{Stealth[length=3mm, width=2mm]}] (6,0)--(8,0);
\draw[-{Stealth[length=3mm, width=2mm]}] (6,0)--(5,{sqrt(3)});
\draw[-{Stealth[length=3mm, width=2mm]}] (6,0)--(5,-{sqrt(3)});
\fill (-6,0) circle (3pt);\fill (0,0) circle (3pt);\fill (6,0) circle (3pt);\node at (-5.1,0.2) {\tiny{$1$}};\node at (-6,1) {\tiny{$q-1$}};\node at (-5.8,-0.8){\tiny{$q^2-q$}};\node at(-1.1,-0.7) {\tiny{$q^2-q$}};\node at (-0.7,0.6) {\tiny{$q$}};\node at (0.9,0.2) {\tiny{$0$}}; \node at (5.4,0.6) {\tiny{$0$}}; \node at (6.9,-0.2) {\tiny{$0$}};\node at (5.4,-0.6) {\tiny{$q^2$}};

\draw [-{Stealth[open, length=3mm, width=2mm]}] (-3,-3.5) -- (-1,-3.5);
\node [draw] at (-2,-3) {in};
\draw [-{Stealth[length=3mm, width=2mm]}] (1,-3.5) -- (3,-3.5);
\node [draw] at (2,-3) {out};
\end{tikzpicture}
\end{center}
\caption{Local transition probability}\label{figure}
\end{figure}
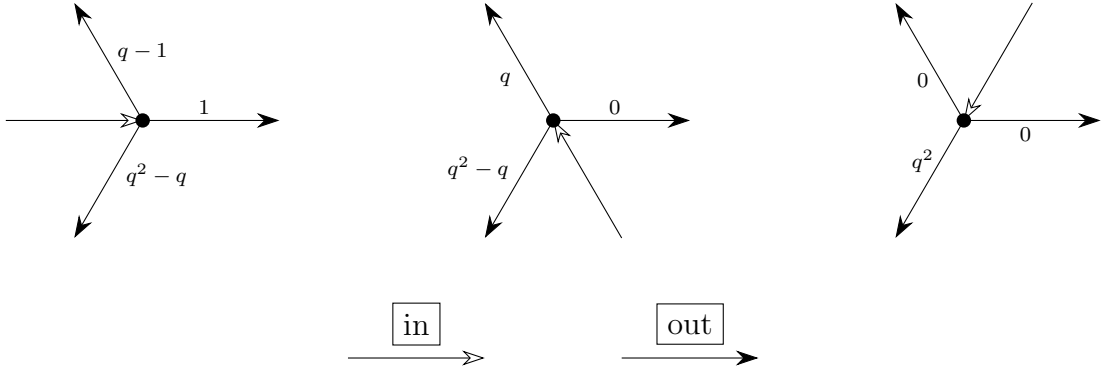

Thus, we obtain the local transition probabilities centered at $v_{m,n}$ of type $1$ geodesic flow on $\Gamma\backslash G/K$ in Figure~\ref{figure}. These probabilities depend on the direction in which they enter toward a fixed vertex. In the case of crossing the boundary of the quotient complex, it will be considered to be going into the reflected vertex against the boundary.

In particular, we have the following lemma.
\begin{lem}
For each positive integer $n$ and a finite segment $(e_0,\ldots,e_{n-1})\in (G/I)^{n}$ of type 1 geodesic in $\mathcal{B}$, 
there are $q^2$ distinct $e_n\in G/I$ such that $(e_0,\ldots,e_{n-1},e_n)\in (G/I)^{n+1}$ is also a type 1 geodesic segment.
\end{lem}
\begin{proof}
This is a special case of Lemma~2.1 in \cite{CM}.
\end{proof}
The directed graph with weights depicted in Figure~\ref{fig:descripD} provides a description of the Markov shift $(\mathcal{D},\sigma)$. In this graph, the (implicitly presented) vertices correspond to the elements of $\Gamma\backslash G/I$, and a directed edge exists from vertex $e_{k,\ell}$ to vertex $e_{k',\ell'}$ only if the pair $(e_{k,\ell}, e_{k',\ell'})$ satisfies $t(e_{k,\ell})=s(e_{k
,\ell'})$. The weights assigned to the edges represent the number of admissible occurrences, that is, a pair $(e_{k,\ell}, e_{k',\ell'})$ lifts to a pair $(e,e')$ of edges in $\mathcal{B}$ such that $s(e),s(e'),t(e')$ do not form a chamber in $\mathcal{B}$.

\begin{center}
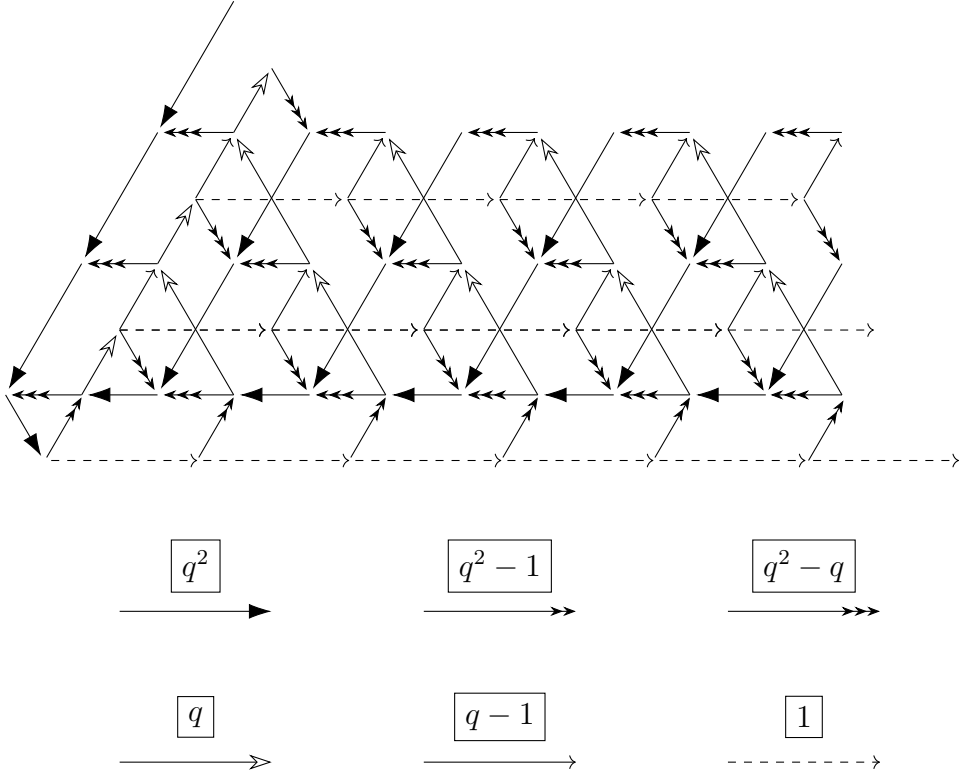
\begin{figure}[h]
\newcommand*\rows{4}
\begin{tikzpicture}[scale=2]
  
\foreach \row in {0, 1, ...,\rows} {
\draw [-{Stealth}{Stealth}] (0.52+\row,0.02) --++ (60:0.46);
\draw [-{Stealth}{Stealth}{Stealth}] (0.749+\row,0.25*1.74) --++ (180:0.46);

\draw [->] (1+\row,0.87) --++ (60:0.46);
\draw [-{Latex[length=3mm, width=2mm]}] (1.249+\row,0.25*1.74) --++ (180:0.46);
\draw [->,dashed] (0.55+\row,0) --++ (0:0.96);
\draw [-{Stealth}{Stealth}{Stealth}] (1.249+\row,0.75*1.74) --++ (180:0.46);
\draw [->] (1.5+\row,1.74) --++ (60:0.46);
\draw [-{Stealth}{Stealth}{Stealth}] (1.749+\row,1.25*1.74) --++ (180:0.46);
\draw [-{Latex[length=3mm, width=2mm]}] (1.25+\row,1.25*1.74) [bend right]--++ (240:0.96);
\draw [-{Latex[length=3mm, width=2mm]}] (0.75+\row,0.75*1.74) --++ (240:0.96);
\draw [-{Stealth[open, length=3mm, width=2mm]}] (1.75+\row,0.25*1.74) --++ (120:0.96);
\draw [->,dashed] (1+\row,0.5*1.73) --++ (0:0.96);

    }
    
\foreach \row in {0, 1, 2, 3}{
\draw [-{Stealth[open, length=3mm, width=2mm]}] (2.25+\row,0.75*1.74) --++ (120:0.96);
\draw [->,dashed] (1.5+\row,1.73) --++ (0:0.96);
\draw [->,dashed] (1+\row,0.5*1.73) --++ (0:0.96);
\draw [-{Stealth}{Stealth}{Stealth}] (1+\row,0.5*1.73) --++ (300:0.46);
\draw [-{Stealth}{Stealth}{Stealth}] (1.5+\row,1.73) --++ (300:0.46);
}
\draw [-{Stealth[open, length=3mm, width=2mm]}] (1.25,0.75*1.74) --++ (60:0.46);
\draw [-{Stealth[open, length=3mm, width=2mm]}] (1.75,1.25*1.74) --++ (60:0.46);
\draw [-{Stealth[open, length=3mm, width=2mm]}] (0.75,0.25*1.74) --++ (60:0.46);
\draw [-{Stealth}{Stealth}] (5.53,0) --++ (60:0.46);
\draw [-{Stealth}{Stealth}{Stealth}] (5.75,0.25*1.74) --++ (180:0.47);
\draw [-{Stealth}{Stealth}{Stealth}] (5,0.5*1.74) --++ (300:0.47);
\draw [-{Latex[length=3mm, width=2mm]}] (0.25,0.25*1.74) --++ (300:0.47);
\draw [-{Stealth}{Stealth}{Stealth}] (5.5,1.73) --++ (300:0.47);
\draw [-{Latex[length=3mm, width=2mm]}] (1.75,1.75*1.74) --++ (240:0.97);
\draw [-{Latex[length=3mm, width=2mm]}] (5.75,0.75*1.74) --++ (240:0.97);
\draw [->,dashed] (5.56,0) --++ (0:0.96);

\draw [-{Latex[length=3mm, width=2mm]}] (1,-1) -- (2,-1);
\node [draw] at (1.5,-0.7) {$q^2$};
\draw [-{Stealth}{Stealth}] (3,-1) -- (4,-1);
\node [draw] at (3.5,-0.7) {$q^2-1$};
\draw [-{Stealth}{Stealth}{Stealth}] (5,-1)--(6,-1); 
\draw [-{Stealth}{Stealth}{Stealth}] (2,2.6)--+(300:0.47); 
\node [draw] at (5.5,-0.7) {$q^2-q$};
\draw [-{Stealth[open, length=3mm, width=2mm]}] (1,-2) -- (2,-2);
\draw [->] (3,-2) -- (4,-2);
\draw [->,dashed] (5,-2) -- (6,-2);
\node [draw] at (1.5,-1.7) {$q$};
\node [draw] at (3.5,-1.7) {$q-1$};
\node [draw] at (5.5,-1.7) {$1$};

\end{tikzpicture}
\caption{Description of $\mathcal{D}$ by directed graph with weights}\label{fig:descripD}
\end{figure}
\end{center}

Let $e_{\frac{1}{2},0}$ be the type $1$ directed edge with $s(e_{\frac{1}{2},0})=v_{0,0}$ and $t(e_{\frac{1}{2},0})=v_{1,0}$. We recall that 
$$\mathcal{D}_{\textrm{per},n}=\{\mathbf{d}\in\mathcal{D}\colon d_0=e_{\frac{1}{2},0},\sigma^n(\mathbf{d})=\mathbf{d}\}$$ 
and
\[g_n(o)=\sum_{\mathbf{d}\in\mathcal{D}_{\textrm{per},n}}\bigr|p_3^{-1}(\mathbf{d})\bigr|.\]
\begin{prop}\label{prop:gn} For each positive integer $n$, we have $$g_{3n}(o)=q^{6n-4}(q^2-1)(q^2-q).$$
\end{prop}
\begin{proof}
We use induction about the distribution of end points in $\mathcal{D}$ for type 1 geodesic segments in $\mathcal{B}$ of length $3n+1$. Let $$N_{\mathcal{D},3n}(e_{k,\ell})=\sum_{\{\mathbf{d}\in\mathcal{D}\colon d_0=e_{\frac{1}{2},0},d_{3n}=e_{k,\ell}\}}\big|p_3^{-1}(d_0,d_1,\ldots,d_{3n})\big|.$$
Then, for each $(k,\ell)$ with $k+\ell\in\frac{1+3\mathbb{Z}}{2}$, we have $$N_{\mathcal{D},3n}(e_{k,\ell})=
\left\{\begin{array}{ll} 
q^{6n-3-2k}(q^2-1)(q^2-q)& \textrm{if }\frac{1}{2}\le k\le 3n-\frac{5}{2},\ell=0\\ 
q^{6n-3-2k}(q^2-1)(q^2-q)& \textrm{if }\frac{5}{2}\le k\le 3n-\frac{5}{2}, \ell=k\\
 q^{6n-2-2k}(q^2-1)^2& \textrm{if }k+\ell<3n-1, 0<\ell<k\\ 
 q^{2\ell}(q^2-1)^2 & \textrm{if }k+\ell=3n-1,\ell\in\frac{1}{2}+\mathbb{Z}, \ell<k \\ 
 q^{2\ell-1}(q^2-1)(q^2-q) & \textrm{if }(k,\ell)=(\frac{3n-1}{2},\frac{3n-1}{2}) \\ 
 q^{2\ell-1}(q^2-1) &\textrm{if }k+\ell=3n+\frac{1}{2},\ell\ne 0 \\ 
1 &\textrm{if }(k,\ell)=(3n+\frac{1}{2},0) \\ 
0 &\textrm{if }k+\ell>3n+\frac{1}{2}
\end{array}\right..$$
Assume that all the equations $N_{\mathcal{D},3n}(e_{k,\ell})$ are correct. Then, it can be readily checked with Figure~\ref{fig:descripD} that the formulas $N_{\mathcal{D},3n+3}(e_{k,\ell})$ are also consistent with the above expressions.
For example, it follows from Figure~\ref{fig:descripD} that
\begin{align*}
N_{\mathcal{D},3n+3}(e_{\frac{1}{2},0})=&\,q^2(q^2-1)(q^2-q)N_{\mathcal{D},3n}(e_{\frac{1}{2},0})+q^4(q^2-q)N_{\mathcal{D},3n}(e_{\frac{3}{2},\frac{1}{2}})\\
&+q^4(q^2-q)N_{\mathcal{D},3n}(e_{2,\frac{3}{2}})+q^6N_{\mathcal{D},3n}(e_{\frac{5}{2},\frac{5}{2}})\\
=&\, q^{6n+2}(q^2-1)(q^2-q).
\end{align*}
Since $g_{3n}(o)=N_{\mathcal{D},3n}(e_{\frac{1}{2},0})$, we get the result.
\end{proof}

\begin{coro}\label{coro:fn} For each positive integer $n$, we have $$f_{3n}(o)=q^{3n-1}(q^2-1)(q^2-q)(q^2+q-1)^{n-1}.$$
\end{coro}
\begin{proof}
This follows directly from the identity $$f_{3n}(o)=g_{3n}(o)-\sum_{k=1}^{n-1}f_{3k}(o)g_{3n-3k}(o)$$
and induction on $n$.
\end{proof}
Proposition~\ref{prop:gn} and Corollary~\ref{coro:fn} yields
\[h_{T_a}=2\log q\qquad\textrm{ and }\qquad \underset{n\to\infty}{\limsup}\,\frac{1}{n}\log f_n(o)=\frac{5}{3}\log q.\] Hence, the inequality in Theorem~\ref{thm:main} directly follows. 

\begin{remark}
The similar directed graph associated to the discrete geodesic flow on $PGL(2,\mathbb{F}_q[t])\backslash PGL(2,\mathbb{F}_q(\!(t^{-1})\!))/PGL(2,\mathbb{F}_q[\![t^{-1}]\!])$ is presented in \cite{Kw}. Defining $g_n(o)$ and $f_n(o)$ similarly, we obtain 
\[g_{2n}(o)=q^{2n-1}(q-1)\quad\textrm{ and }\quad f_{2n}(o)=q^n(q-1).\] In particular, $f_{2n}(o)$ is equal to the number of degree $n$ polynomials in $\mathbb{F}_q[t]$. These values are related to the partial quotients of the continued fraction expansion of quadratic irrationals in $\mathbb{F}_q(\!(t^{-1})\!)$ (cf. \cite{PS}). We believe that it would be interesting to discover an alternative interpretation of the formula of $f_{3n}(o)$ in $PGL_3$ through multi-dimensional continued fraction theory.
\end{remark}

\end{document}